\documentclass{gtart}

\def\ifplaintex{\expandafter\ifx\csname documentclass\endcsname\relax}

\def\gtp{{\mathsurround=0pt\it $\cal G\mskip-2mu$eometry \&\ 
$\cal T\!\!$opology $\cal P\!$ublications}}  

\def\recd{{\small Received:\qua\receiveddate\ifx\reviseddate\relax
\else\qquad Revised:\qua\reviseddate\fi\par}} 

\def\lognumber#1{\def\thelognumber{#1}}
\def\volumenumber#1{\def\thevolumenumber{#1}}
\def\volumeyear#1{\def\thevolumeyear{#1}}
\def\papernumber#1{\def\thepapernumber{#1}}
\def\pagenumbers#1#2{\def\startpage{#1}\def\finishpage{#2}}
\def\published#1{\def\publishdate{#1}}

\def\received#1{\def\receiveddate{#1}}

\def\accepted#1{\def\accepteddate{#1}}
\def\asciititle#1{\def\theasciititle{#1}}

\def\asciiaddress#1{\def\theasciiaddress{#1}}

\long\def\asciiabstract#1{\long\def\theasciiabstract{#1}}

\let\\\par\let\thelognumber\relax\let\thevolumenumber\relax
\let\thepapernumber\relax\let\thevolumeyear\relax\let\startpage\relax
\let\finishpage\relax\let\publishdate\relax\let\receiveddate\relax
\let\reviseddate\relax\let\accepteddate\relax\let\theasciititle\relax
\let\theasciiauthors\relax\let\theasciiaddress\relax
\let\theasciiabstract\relax

\let\theasciiemail\relax

\ifplaintex
\font\logobig=cmssbx10 scaled 3836
\font\logomed=cmssbx10 scaled 2557
\else
\font\logobig=cmssbx10 scaled 4200
\font\logomed=cmssbx10 scaled 2800
\fi

\long\def\makeagttitle{   
\count0=\startpage
\agt\hfill      
\hbox to 45truept{\vbox to 0pt{\vglue -13truept{\logomed A\kern -.37em{\logobig 
T}\kern -.38em G}\vss}\hss}
\break
{\small Volume \thevolumenumber\ (\thevolumeyear)
\startpage--\finishpage\nl
Published: \publishdate}

\vglue .25truein

{\parskip=0pt\leftskip 0pt plus
1fil\def\\{\par\smallskip}{\Large\bf\thetitle}\par\medskip} \vglue
0.05truein

{\parskip=0pt\leftskip 0pt plus 1fil\def\\{\par}{\sc\theauthors}
\par\medskip}%
 
\vglue 0.03truein 

{\small\leftskip 25truept\rightskip 25truept{\bf Abstract}\stdspace\theabstract

{\bf AMS Classification}\stdspace\theprimaryclass
\ifx\thesecondaryclass\relax\else; \thesecondaryclass\fi\par
{\bf Keywords}\stdspace \thekeywords\par}\vglue 7truept

}   

\ifplaintex
\font\phead=cmsl9 scaled 950
\font\pnum=cmbx10 scaled 913
\font\pfoot=cmsl9 scaled 950
\headline{\vbox to 0pt{\vskip -4.5mm\line{\small\phead\ifnum
\count0=\startpage ISSN 1472-2739 (on-line) 1472-2747 (printed)
\hfill {\pnum\folio}\else\ifodd\count0\def\\{ }%
\ifx\theshorttitle\relax\thetitle\else\theshorttitle\fi\hfill{\pnum\folio}
\else\def\\{ and }{\pnum\folio}\hfill\ifx\theshortauthors\relax\theauthors
\else\theshortauthors\fi\fi\fi}\vss}}
\footline{\vbox to 0pt{\vglue 0mm\line{\small\pfoot\ifnum\count0=\startpage
\copyright\ \gtp\hfill\else
\agt, Volume \thevolumenumber\ (\thevolumeyear)\hfill\fi}\vss}}
\else
\font=cmsl9 scaled 1050
\font\lnum=cmbx10 
\font=cmsl9 scaled 1050
\makeatletter
\def\@oddhead{{\small\lhead\ifnum\count0=\startpage ISSN 1472-2739 
(on-line) 1472-2747 (printed)\hfill {\lnum\number\count0}\else\ifodd\count0
\def\\{ }\ifx\theshorttitle\relax \thetitle \else\theshorttitle\fi\hfill
{\lnum\number\count0}\else\def\\{ and }{\lnum\number\count0}
\hfill\ifx\theshortauthors\relax 
\theauthors\else\theshortauthors\fi\fi\fi}}\def\@evenhead{\@oddhead}
\def\@oddfoot{\small\lfoot\ifnum\count0=\startpage\copyright\ \gtp\hfill\else
\agt, Volume \thevolumenumber\ (\thevolumeyear)\hfill\fi}
\def\@evenfoot{\@oddfoot}
\makeatother
\fi
\let\maketitlepage\makeagttitle
\let\makeshorttitle\maketitlepage
\let\maketitle\maketitlepage

\newwrite\gtoutfile
\long\gdef\makeheadfile{  
{\def\\{, }\def\s{ }
\immediate\openout\gtoutfile head.xxx
\immediate\write\gtoutfile{Proxy-for: \ifx\theasciiauthors\relax
\theauthors\else\theasciiauthors\fi\s<\ifx\theasciiemail\relax\theemail\else\theasciiemail\fi>}
\immediate\write\gtoutfile{\noexpand\\}
\immediate\write\gtoutfile{Authors: \ifx\theasciiauthors\relax
\theauthors\else\theasciiauthors\fi}
{\def\\{ }\immediate\write\gtoutfile{Title: \ifx\theasciititle\relax
\thetitle\else\theasciititle\fi}}
\immediate\write\gtoutfile{Subj-class: GT or SG, GR etc}
\immediate\write\gtoutfile{MSC-class: \theprimaryclass\ifx\thesecondaryclass\relax\else, \thesecondaryclass\fi}
\immediate\write\gtoutfile{Journal-ref: Algebr. Geom. Topol. \thevolumenumber\s
(\thevolumeyear) \startpage-\finishpage}
\immediate\write\gtoutfile{Comments: Published by Algebraic and
Geometric Topology at}
\immediate\write\gtoutfile{\s\s\s  http://www.maths.warwick.ac.uk/agt/AGTVol\thevolumenumber/agt-\thevolumenumber-\thepapernumber.abs.html}
\immediate\write\gtoutfile{\noexpand\\}
\immediate\write\gtoutfile{}
\ifx\theasciiabstract\relax
\immediate\write\gtoutfile{\theabstract}\else
\immediate\write\gtoutfile{\theasciiabstract}\fi
\immediate\write\gtoutfile{}
\immediate\write\gtoutfile{\noexpand\\}
\immediate\write\gtoutfile{}
\immediate\closeout\gtoutfile}}  

\def\maketitlepage{\makeagttitle\makeheadfile}
\let\makeshorttitle\maketitlepage
\let\maketitle\maketitlepage

\lognumber{19}
\volumenumber{3}
\volumeyear{2003}
\papernumber{19}
\published{20 June 2003}
\pagenumbers{569}{586}
\received{27 January 2003}
\accepted{23 April 2003}

\usepackage{amssymb,graphicx,amsmath}

\theoremstyle{definition}
\newtheorem{defn}{Definition}[section]

\newtheorem{rmk}[defn]{Remark}

\newtheorem{qn}[defn]{Question}
\theoremstyle{plain}
\newtheorem{thm}[defn]{Theorem} 
\newtheorem{cor}[defn]{Corollary}

\newtheorem{correction}{Correction}

\newcommand{\reals}{\ensuremath{\mathbb{R}}}
\newcommand{\half}{\ensuremath{\frac{1}{2}}}
\newcommand{\into}{\ensuremath{\hookrightarrow}}

\newcommand{\pf}{\mathop{\rm pf}\nolimits}
\newcommand{\cf}{\mathop{\rm cf}\nolimits}

\newcommand{\genus}{\mathop{\rm genus}\nolimits}
\begin{document}
\title{Open books and configurations\\of symplectic surfaces}
\asciititle{Open books and configurations of symplectic surfaces and erratum}
\authors{David T. Gay}                  
\address{Department of Mathematics, University of Arizona\\
617 North Santa Rita, PO Box 210089\\
Tucson, AZ 85721, USA}                  
\email{dtgay@math.arizona.edu}

\begin{abstract} 
We study neighborhoods of configurations of symplectic surfaces in
symplectic $4$--manifolds. We show that suitably ``positive''
configurations have neighborhoods with concave boundaries and we
explicitly describe open book decompositions of the boundaries
supporting the associated negative contact structures. This is used to
prove symplectic nonfillability for certain contact $3$--manifolds and
thus nonpositivity for certain mapping classes on surfaces with
boundary. Similarly, we show that certain pairs of contact
$3$--manifolds cannot appear as the disconnected convex boundary of
any connected symplectic $4$--manifold. Our result also has the
potential to produce obstructions to embedding specific symplectic
configurations in closed symplectic $4$--manifolds and to generate new
symplectic surgeries. From a purely topological perspective, the
techniques in this paper show how to construct a natural open book
decomposition on the boundary of any plumbed $4$--manifold.
\end{abstract}

\asciiabstract{We study neighborhoods of configurations of symplectic
surfaces in symplectic 4-manifolds. We show that suitably
`positive' configurations have neighborhoods with concave boundaries
and we explicitly describe open book decompositions of the boundaries
supporting the associated negative contact structures. This is used to
prove symplectic nonfillability for certain contact 3-manifolds and
thus nonpositivity for certain mapping classes on surfaces with
boundary. Similarly, we show that certain pairs of contact
3-manifolds cannot appear as the disconnected convex boundary of
any connected symplectic 4-manifold. Our result also has the
potential to produce obstructions to embedding specific symplectic
configurations in closed symplectic 4-manifolds and to generate new
symplectic surgeries. From a purely topological perspective, the
techniques in this paper show how to construct a natural open book
decomposition on the boundary of any plumbed 4-manifold.

Erratum (added December 2003): 
We correct the main theorem and its proof.  As originally stated, the
theorem gave conditions on a configuration of symplectic surfaces in a
symplectic 4-manifold under which we could construct a model
neighborhood with concave boundary and describe explicitly the open
book supporting the contact structure on the boundary. The statement
should have included constraints on the areas of the surfaces.}

\primaryclass{57R17}                
\secondaryclass{57N10, 57N13}              

\maketitle

\section{Main Results}
All manifolds in this paper are oriented; if $(X,\omega)$ is a
symplectic $4$--manifold we assume that $\omega \wedge \omega > 0$.
By a {\em symplectic configuration} in a symplectic $4$--manifold
$(X,\omega)$, we mean a union $C = \Sigma_1 \cup \ldots \cup \Sigma_n$
of closed symplectic surfaces embedded in $(X,\omega)$ such that all
intersections between surfaces are $\omega$--orthogonal.

A {\em symplectic configuration graph} is a labelled graph $G$ with no
edges from a vertex to itself and with each vertex $v_i$ labelled with
a triple $(g_i,m_i,a_i)$, where $g_i \in \{0,1,2,\ldots\}$, $m_i \in
\mathbb{Z}$ and $a_i \in (0,\infty)$.  Associated to a symplectic
configuration $C = \Sigma_1 \cup \ldots \cup \Sigma_n$ in a symplectic
$4$--manifold $(X,\omega)$ is a symplectic configuration graph $G(C)$
where each vertex $v_i$ corresponds to a surface $\Sigma_i$, $g_i =
\genus(\Sigma_i)$, $m_i = \Sigma_i \cdot \Sigma_i$ and $a_i =
\int_{\Sigma_i} \omega$, and where each edge represents a point of
intersection. Because $\omega$--orthogonal intersections are
necessarily positive, $G(C)$ completely determines the topology of a
regular neighborhood of $C$ (namely, the result of {\em plumbing} disk
bundles over surfaces according to $G(C)$); for this much the areas
$\{a_i\}$ are irrelevant.  If we include the area information then
$G(C)$ also determines the germ of $\omega$ near $C$ (due to a slight
generalization of standard symplectic neighborhood theorems, as
in~\cite{Symington}).

For any vertex $v_i$ in a graph $G$, let $d_i$ denote the degree of
$v_i$, the number of edges connected to $v_i$. We say that a
configuration graph $G$ is {\em positive} if $m_i + d_i > 0$ for every
vertex $v_i$.

Recall that the boundary of a symplectic $4$--manifold $(X,\omega)$ is
{\em concave} (resp.\ {\em convex}) if there exists a symplectic
dilation $V$ defined on a neighborhood of $\partial X$ pointing in
(resp.\ out) along $\partial X$; this induces a negative
(resp.\ positive) contact structure $\xi = \ker \imath_V \omega
|_{\partial X}$ on $\partial X$.

We present our main theorem in two parts. Part~A states that positive
symplectic configurations have neighborhoods with concave boundaries,
and part~B explicitly describes the contact structures on such
boundaries in terms of open book decompositions.

\begin{thm}{\bf(part A)}\label{MainThm}\ \
Given any positive symplectic configuration graph~$G$ there exists an open
symplectic $4$--manifold $(N(G),\omega(G))$, a symplectic
configuration $C(G) \subset (N(G),\omega(G))$ and a positive contact
$3$--manifold $(M(G),\xi(G))$, satisfying the following properties:
\begin{itemize}
\item $G=G(C(G))$.
\item For some contact form $\alpha$ for $\xi(G)$ and some (not
necessarily smooth) function $f: M(G) \to \reals$, letting $X_f =
\{(t,p) | t < f(p)\} \subset \reals \times M(G)$ and $\omega = d(e^t
\alpha)$, there is a symplectomorphism $\phi: (X_f,\omega) \to (N(G)
\setminus C(G),\omega(G))$ such that
\[ C(G) = \{ \lim_{t \rightarrow f(p)} \phi(t,p) | p \in M(G) \}. \]
\end{itemize} 
Thus, given any symplectic configuration $C$ in any symplectic
$4$--manifold $(X,\omega)$, if $G(C) = G$ then there exists a compact
neighborhood of $C$ in $(X,\omega)$ which is symplectomorphic to a
neighborhood of $C(G)$ in $(N(G),\omega(G))$ and
which has concave boundary contactomorphic to $(-M(G),\xi(G))$.
\end{thm}

The strength of this theorem will lie in the characterization of
$(M(G),\xi(G))$ in terms of an open book decomposition of $M(G)$. We
briefly recall the relationship between contact structures and open
books; for more details see~\cite{Giroux} and~\cite{GayCF}.

An open book decomposition of a $3$--manifold $M$ is a pair $(L,p)$,
where $L$ is a link and $p: M \setminus L \rightarrow S^1$ is a
fibration such that the fibers are longitudinal near each component of
$L$. The link $L$ is called the {\em binding} and the compact surfaces
$\Sigma_t = p^{-1}(t) \cup L$ are called the {\em pages}, with $L =
\partial \Sigma_t$ for all $t \in S^1$. By the mapping class group
$\mathcal{M}(\Sigma)$ for a compact surface $\Sigma$ with boundary, we
mean the group of orientation-preserving self-diffeomorphisms of
$\Sigma$ fixing $\partial \Sigma$ pointwise modulo isotopies fixing
$\partial \Sigma$ pointwise. The {\em monodromy} of an open book is
the mapping class $h \in \mathcal{M}(\Sigma_0)$ given by the return
map of a flow transverse to the pages and meridinal near the binding.

A positive contact form $\alpha$ on $M$ is {\em supported} by the open
book $(L,p)$ if $d\alpha$ is positive on each page and if $\alpha$
orients $L$ in the same sense that $L$ is oriented as the boundary of
a page. A positive contact structure $\xi$ is {\em supported} by
$(L,p)$ if $\xi = \ker \alpha$ for some contact form $\alpha$ which is
supported by $(L,p)$. We have the following result at our disposal:
\begin{thm}[Thurston-Winkelnkemper~\cite{ThurstonWinkel}, 
Torisu~\cite{Torisu}, Giroux~\cite{Giroux}] \label{TWTG}
Every open book decomposition of any $3$--manifold supports some
positive contact structure, and any two positive contact structures
supported by the same open book are isotopic.
\end{thm}
Thus, given a compact surface $\Sigma$ with boundary and a mapping
class $h \in \mathcal{M}(\Sigma)$, there exists a unique (up to
contactomorphism) positive contact $3$--manifold with contact structure
supported by an open book with page $\Sigma$ and monodromy $h$; we
denote this contact manifold $\mathcal{B}(\Sigma,h)$.

\begin{figure}
\begin{center}
\includegraphics[width=4in,height=1.5in]{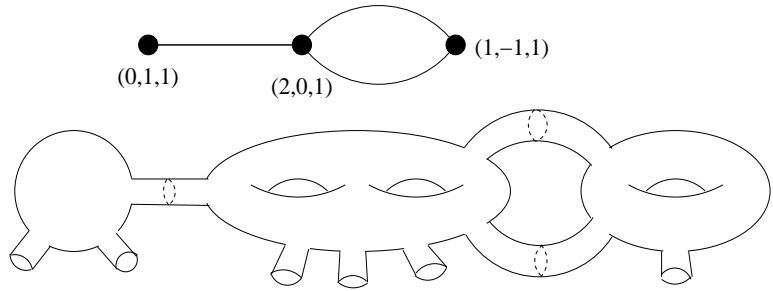}
\caption{A symplectic configuration graph $G$ and the surface $\Sigma(G)$}
\label{F:SigmaGExample}
\end{center}
\end{figure}

Given a positive configuration graph $G$, for each vertex $v_i$ let
$F_i$ be a surface of genus $g_i$ with $m_i + d_i$ boundary
components. Let $\Sigma(G)$ be the surface resulting from performing
connect sums between these surfaces, with one connect sum between
$F_i$ and $F_j$ for each edge connecting $v_i$ to $v_j$. Each edge in
$G$ corresponds to a circle in $\Sigma(G)$. An example of a graph $G$
and the surface $\Sigma(G)$ is illustrated in
figure~\ref{F:SigmaGExample}, with the circles corresponding to the
edges drawn in dashed lines.
Let $\sigma(G)$ be the product of one right-handed Dehn twist around
each of the circles in $\Sigma(G)$ corresponding to the edges of $G$,
let $\delta(G)$ be the product of one right-handed Dehn twist around
each component of $\partial \Sigma(G)$ and let $h(G) = \sigma(G)^{-1}
\circ \delta(G)$.

\vspace{6pt}
\begin{flushleft}
\textbf{Theorem~\ref{MainThm}\qua (part B)}\qua
$(M(G), \xi(G)) = \mathcal{B}(\Sigma(G),h(G))$
\end{flushleft}
\vspace{6pt}

Note that the area information recorded in the graph $G$ is not
recorded on the boundary $(M(G),\xi(G))$. For this reason, we may
suppress mention of the areas and think of the vertices of $G$ as
labelled with pairs $(g,m)$, rather than triples $(g,m,a)$.

After proving theorem~\ref{MainThm} we will briefly discuss the
situation when $G$ is not positive. 

Theorem~\ref{MainThm} also has purely topological content, namely the
explicit characterization of an open book on the boundary of a plumbed
$4$--manifold corresponding to a positive configuration graph. If we
are only interested in smooth topology, the techniques used to prove
theorem~\ref{MainThm} do yield a theorem characterizing open books on
the boundaries of arbitrary plumbed $4$--manifolds, which we now
state.

Let a {\em plumbing graph} be a graph $G$ with no edges connecting a
vertex to itself, with each vertex $v_i$ labelled with a pair
$(g_i,m_i)$ and with each edge labelled with a $+$ or a $-$. The
plumbed $4$--manifold $X(G)$ corresponding to $G$ is a neighborhood of
a configuration of surfaces $\Sigma_1 \cup \ldots \cup \Sigma_n$
corresponding to the vertices $v_1, \ldots, v_n$ of $G$, with
$\genus(\Sigma_i) = v_i$, $\Sigma_i \cdot \Sigma_i = m_i$ and each $+$
(resp.\ $-$) edge corresponding to a positive (resp.\ negative)
transverse intersection between two surfaces. For each vertex $v_i$
let $d_i$ be the {\em signed} count of edges connecting to $v_i$ (a
$+$ edge contributes $+1$ while a $-$ edge contributes $-1$). For each
$v_i$ let $F_i$ be a surface of genus $g_i$ with $|m_i + d_i|$
boundary components and let $\Sigma(G)$ be the result of performing
connect sums between these surfaces according to $G$. Let $h(G)$ be
the product of the following Dehn twists: one right-handed Dehn twist
about each circle in $\Sigma(G)$ corresponding to a $+$ edge in $G$,
one left-handed Dehn twist about each circle in $\Sigma(G)$
corresponding to a $-$ edge, one left-handed Dehn twist about each
boundary component coming from a vertex $v_i$ for which $m_i + d_i >
0$, and one right-handed Dehn twist about each boundary component
coming from a vertex for which $m_i + d_i < 0$.

\begin{thm} \label{TopThm}
Given any plumbing graph $G$, let $X(G)$ be the associated plumbed
$4$--manifold. Then $\partial X$ has an open book decomposition with
page $\Sigma(G)$ and monodromy $h(G)$.
\end{thm}
In the case where $m_i + d_i = 0$ for all vertices, this is in fact
giving us a fibration of $\partial X$ over $S^1$, or an ``open book
with empty binding''. If we apply theorem~\ref{TopThm} to a positive
symplectic configuration graph, the reader may notice that the
monodromy as described here is the inverse of the monodromy as
described in theorem~\ref{MainThm}; this is because here we are
describing $\partial X$ whereas in theorem~\ref{MainThm} we are
describing $-\partial X$.

The author would like to thank A. Stipsicz for suggesting the idea of
trying to understand boundary behaviors for neighborhoods of
symplectic configurations as a way to search for new symplectic
surgeries, and would like to thank A. Stipsicz, G. Matic, M. Symington
and R. Kirby for helpful discussions and for looking at drafts of this
paper and suggesting improvements.

\section{Applications}

Before presenting the main proofs we investigate a few consequences of
theorem~\ref{MainThm} and point out some directions in which to look
for further applications.

Given a compact surface $\Sigma$, we say that a mapping class $h \in
\mathcal{M}(\Sigma)$ is positive if $h$ can be expressed as a product
of right-handed Dehn twists. It is not hard to show, using compact
Stein surfaces and Legendrian surgeries, that if $h \in
\mathcal{M}(\Sigma)$ is positive then $\mathcal{B}(\Sigma,h)$ is
strongly symplectically fillable (see~\cite{LoiPier}, \cite{AkbOzb}
and~\cite{GayCF}).

Given a configuration graph $G$, let $Q(G)$ be the associated
intersection form; i.e.\ $Q(G) = (q_{ij})$, where $q_{ii} = m_i$ and
$q_{ij}$ is the number of edges connecting $v_i$ to $v_j$. Let
$b^+(G)$ denote the number of positive eigenvalues of $Q(G)$.

The following is a straightforward application of the adjunction
inequality (see~\cite{McDSal}):
\begin{cor} \label{C1}
Let $G$ be a connected positive graph with $b^+(G) > 1$ and with at
least one vertex $v_i$ for which $m_i > 2g_i - 2$. Then
$\mathcal{B}(\Sigma(G),h(G))$ is not strongly symplectically fillable
and therefore $h(G)$ is not positive in $\mathcal{M}(\Sigma(G))$.
\end{cor}

\begin{proof}[Proof of corollary~\ref{C1}]
The symplectic manifold $(N(G),\omega(G))$ constructed in
theorem~\ref{MainThm} is open; a function $F < f$ on $M(G)$ gives a
compact version $(N_F(G),\omega(G))$ where $N_F(G) = \phi\{(t,p) |
F(p) \leq t < f(p)\} \cup C(G)$.  Suppose that $\mathcal{B}(\Sigma(G),
h(G))$ is strongly symplectically fillable. Then there exists a closed
symplectic $4$--manifold $(X,\omega)$ containing $(N_F(G),\omega(G))$
for some function $F < f$ on $M(G)$ and containing the configuration
$C(G) = \Sigma_1 \cup \ldots \cup \Sigma_n$. The intersection form for
$N_F(G)$ is $Q(G)$; since $b^+(G) > 1$ we know that $b_2^+(X) >
1$. Thus the adjunction inequality applies, which states that, for any
closed surface $\Sigma \subset X$ (with $[\Sigma]$ not torsion in
$H_2(X)$ if $\genus(\Sigma) = 0$), $|c_1(\omega) \cdot \Sigma| +
\Sigma \cdot \Sigma \leq 2 \genus(\Sigma) - 2$. However, here we have
an embedded surface $\Sigma_i$ for which $\Sigma_i \cdot \Sigma_i >
2\genus(\Sigma_i) - 2$, which is a contradiction.
\end{proof}

\begin{rmk} \label{R1}
For any surface $\Sigma$, let us call a relation in $\mathcal{M}(G)$
of the form $\delta = w$ a {\em boundary-interior relation} if
$\delta$ is a single right twist about each boundary component and $w$
is some word in interior right twists.  We have the following trivial
observation: For a given $G$, $h(G)$ is positive if and only if there
exists a boundary-interior relation $\delta(G) = w$ in
$\mathcal{M}(\Sigma(G))$ such that the word $w$ includes all the
twists in $\sigma(G)$. (The order in which the twists of $\sigma(G)$
appear in $w$ does not matter.) Boundary-interior relations have a
variety of uses, including giving constructions of topological
Lefschetz pencils (see~\cite{GayCF}). 

Let $\Sigma_g^n$ denote a surface of genus $g$ with $n$ boundary
components.  Two boundary-interior relations are the ``lantern
relation'' in $\mathcal{M}(\Sigma_0^4)$ and the ``chain relation'' in
$\mathcal{M}(\Sigma_g^2)$ (see~\cite{Wajnryb}).
\begin{figure}
\begin{center}
\includegraphics[width=4in,height=3in]{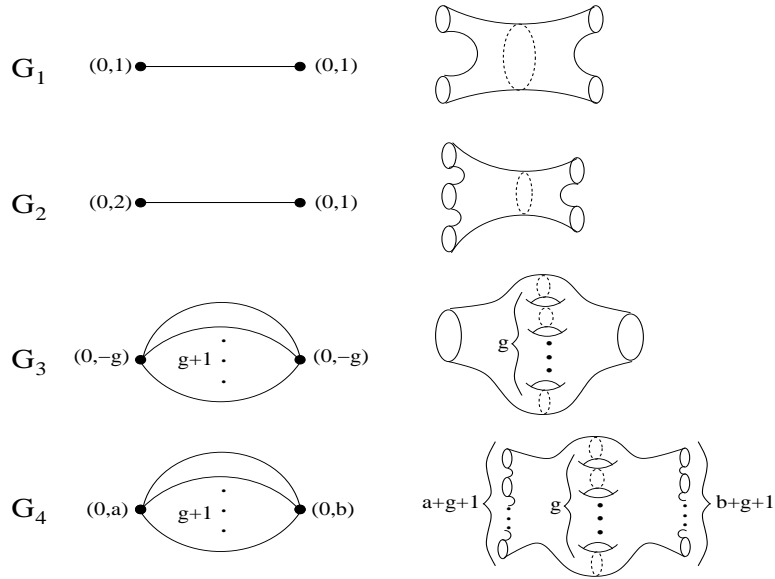}
\caption{Four example graphs for remark~\ref{R1}}
\label{F:RelationExamples}
\end{center}
\end{figure}
Figure~\ref{F:RelationExamples} shows a set of example graphs $G_1$,
$G_2$, $G_3$, $G_4$ on the left (here we have suppressed the areas and
only given the pair $(g,m)$ at each vertex), with the associated
surfaces $\Sigma(G_1), \ldots, \Sigma(G_4)$ drawn on the right. The
lantern relation shows that $h(G_1)$ is positive, while
corollary~\ref{C1} tells us that $h(G_2)$ is not positive. Thus there
does not exist a lantern-type relation on $\Sigma_0^5 =
\Sigma(G_2)$. The chain relation shows that $h(G_3)$ is positive
(where $G_3$ has $g+1$ edges between $2$ vertices so that $\Sigma(G_3)
= \Sigma_g^2$). The graph $G_4$ also has $g+1$ edges connecting two
vertices, but now the self-intersections are arbitrary integers $a$
and $b$ such that $a + g + 1 > 0$ and $b+g+1 > 0$. If $a$ and $b$ are
positive and $ab > (g+1)^2$, then $b^+(G_4) = 2$ and
corollary~\ref{C1} tells us that $h(G_4)$ is not positive. Thus, if
$ab > (g+1)^2$, $A = a + g + 1$, $B = b+g+1$ and $n = A+B$, then there
cannot exist a boundary-interior relation $\delta = w$ in
$\mathcal{M}(\Sigma_g^n)$ if $w$ contains twists along $g+1$ disjoint
curves which, collectively, separate $\Sigma_g^n$ into two genus $0$
pieces one containing $A$ of components of $\partial \Sigma_g^n$ and
the other containing $B$ components.

The existence of an elliptic Lefschetz pencil on $\mathbb{CP}^2$ with
$9$ points in the base locus and $12$ singular fibers means that there
exists a boundary-interior relation $\delta = w$ in
$\mathcal{M}(\Sigma_1^9)$ such that $w$ is the product of twists along
$12$ curves $C_1, \ldots, C_{12}$. This author is not aware that this
relation has been written down and has been curious for a long time as
to where these $12$ curves are.  We do know that when we blow up the
$9$ points we get a Lefschetz fibration given by the relation $(ab)^6
= 1$ in $\mathcal{M}(\Sigma_1^0)$, where $a$ is a meridinal right
twist and $b$ is a longitudinal right twist, so that, after embedding
$\Sigma_1^9$ in $\Sigma_1^0$ the odd $C_i$'s are isotopic to meridians
and the even $C_i$'s are isotopic to longitudes. Corollary~\ref{C1}
gives a little more information, ruling out certain possible
arrangements of curves. For example, the observation from the previous
paragraph about $G_4$ shows that no two of the curves may be disjoint
and separate $4$ boundary components from the other $5$. Other
possibilities can be ruled out by considering various cycle graphs.
\end{rmk}

Along much the same lines, we have:
\begin{cor}
Let $G_1$ and $G_2$ be positive, connected symplectic configuration
graphs with $b^+(G_1) > 0$ and $b^+(G_2) > 0$, with at least one
vertex $v_i$ in $G_1$ for which $m_i > 2g_i - 2$. Then, although each
$\mathcal{B}(\Sigma(G_i),h(G_i))$ may be strongly symplectically
fillable, there does not exist a connected symplectic $4$--manifold
with disconnected convex boundary $\mathcal{B}(\Sigma(G_1),h(G_1))
\amalg \mathcal{B}(\Sigma(G_2),h(G_2))$.
\end{cor}
\begin{proof}
If such a symplectic $4$--manifold existed then there would exist a
closed, connected, symplectic $4$--manifold $(X,\eta)$ containing
\[(N,\omega) = (N_{F_1}(G_1),\omega(G_1)) \amalg
(N_{F_2}(G_2),\omega(G_2))\] 
for appropriate functions $F_1$ and
$F_2$. The intersection form for $N$ is $Q(G_1) \oplus Q(G_2)$, so
that $b_2^+(N) > 1$; the rest of the contradiction is identical to
that in the preceding proof.
\end{proof}

\begin{rmk}
McDuff~\cite{McDuffCvx} has shown that symplectic $4$--manifolds with
disconnected convex boundary do exist. By a result of
Gromov~\cite{Gromov} (made explicit in~\cite{EliashFilling}
and~\cite{McDuffCvx}), it is not possible to have multiple convex
$S^3$ boundary components, which is the case of this corollary when
$\Sigma(G_1)$ and $\Sigma(G_2)$ are both disks. We hope that our
result significantly enlarges the class of pairs of contact manifolds
which cannot appear as disconnected convex boundaries, despite being
individually strongly symplectically fillable. It is not hard to
construct examples of graphs $G$ with $b^+(G) = 1$, with one vertex
for which $m_i > 2g_i - 2$ and such that $M(G)$ is not $S^3$, but it
is then not necessarily clear how to show that $(M(G),\xi(G))$ is in
fact strongly symplectically fillable.
\end{rmk}

It would be interesting to construct arguments in the opposite
direction:
\begin{qn} Are there any positive configuration graphs $G$ for which
we can show directly that $\mathcal{B}(\Sigma(G),h(G))$ is overtwisted
and hence conclude, without appealing to the adjunction inequality,
that a symplectic configuration with configuration graph $G$ cannot
embed in a closed symplectic $4$--manifold?
\end{qn}
Showing that $\mathcal{B}(\Sigma,h)$ is overtwisted for a given
surface $\Sigma$ and mapping class $h \in \mathcal{M}(\Sigma)$ is, in
principle, purely a mapping class group problem, as follows: Let
$\Sigma$ and $\Sigma'$ be compact surfaces with boundary and $h \in
\mathcal{M}(\Sigma)$ and $h' \in \mathcal{M}(\Sigma')$ be mapping
classes. We say that $(\Sigma',h')$ is a positive (resp.\  negative)
stabilization of $(\Sigma,h)$ if $\Sigma'$ is the result of attaching
a $1$--handle to $\Sigma$ and $h' = h \circ \tau$, where $\tau$ is a
right-handed (resp.\ left-handed) Dehn twist along a curve dual to the
co-core of the $1$--handle.  It can be shown, using results of
Giroux~\cite{Giroux} and Torisu~\cite{Torisu}, that
$\mathcal{B}(\Sigma,h)$ is overtwisted if and only if there exist
pairs $(\Sigma',h')$ and $(\Sigma'',h'')$ such that $(\Sigma',h')$ is
related to $(\Sigma,h)$ by a sequence of positive stabilizations and
destabilizations and $(\Sigma',h')$ is a negative stabilization of
$(\Sigma'',h'')$

Lastly, we point out that theorem~\ref{MainThm} could be used
to produce new symplectic surgeries. If, for a given symplectic
configuration graph $G$, we can find some other interesting symplectic
$4$--manifold $(Y,\eta)$ with concave boundary contactomorphic to
$(-M(G),\xi(G))$, then we may surger out a neighborhood of a
configuration $C$ for which $G(C) = G$ and replace it with
$(Y,\eta)$. (Symington~\cite{Symington, SymingtonGen} has investigated
configurations which have neighborhoods with convex boundaries, and
has used this to produce useful symplectic surgeries; this paper is
partly inspired by that work.)
\begin{qn}
Given $G$, is there any canonical way to produce such a $(Y,\eta)$
with significantly different topology from $(N(G),\omega(G))$?
Especially interesting would be examples where $Y$ is a rational
homology ball.
\end{qn}

\section{The main proof}

\begin{proof}[Proof of theorem~\ref{MainThm}]
Our proof is a three step construction. In ``Step~1'' we construct a
symplectic $4$--manifold $(X,\omega)$ with a symplectic dilation
(Liouville vector field) $V^+$ defined on all of $X$ and pointing out
along $\partial X$. $X$ will contain a configuration $Z$ of properly
embedded symplectic surfaces with boundary, which will become the
configuration $C(G)$ once we cap off the surfaces with $2$--handles. In
``Step~2'', we arrange that the induced positive contact form on
$\partial X$ has a particularly nice form and then we produce a
symplectic contraction $V^-$ defined on $X \setminus Z$, pointing out
along $\partial X \setminus \partial Z$.  In ``Step~3'' we cap off the
surfaces in $Z$ by attaching appropriately framed $4$--dimensional
symplectic $2$--handles along $\partial Z$; these handles have the
effect of turning the convex boundary into a concave boundary. The
symplectic contraction $V^-$ from Step~2 extends across the
$2$--handles; flow along $-V^-$ gives the symplectomorphism $\phi$
(after attaching an open collar to the boundary).

{\bf Step~1}\qua For lack of a better term, the objects we construct in
this step will be called ``Step~1 objects''. A Step~1 object is a
$6$--tuple $(X,\omega,Z,V^+,f,p)$ where:
\begin{itemize}
\item $(X,\omega)$ is a compact symplectic $4$--manifold with boundary.
\item $Z = F_1 \cup \ldots \cup F_n$ is a configuration of symplectic
surfaces with boundary, each properly embedded in $X$, with
$\omega$--orthogonal intersections.
\item $f$ is a proper Morse function on $X$ which restricts to each
  $F_i$ as a proper Morse function, with only critical points of index
  $0$ and $1$, all of which lie in $Z$.
\item $V^+$ is a symplectic dilation on $(X,\omega)$, tangent to $Z$
  and gradient-like for $f$, inducing a positive contact
  structure $\xi^+ = \ker (\imath_{(V^+)} \omega|_{\partial X})$ on
  $\partial X$.
\item $p: \partial X \setminus \partial Z \rightarrow S^1$ is a
fibration such that $(\partial Z,p)$ is an open book on $\partial X$.
\item $\xi^+$ is supported by $(\partial Z,p)$.
\end{itemize}

Each component $K$ of $\partial F_i \subset \partial Z$ has two
different natural framings, the framing coming from the page of the
open book, which is a Seifert surface for $\partial Z$, and the
framing coming from a Seifert surface for $\partial F_i$. Denote
the former framing $\pf(K)$ (for ``page framing'') and denote the
latter framing $\cf(K)$ (for ``component framing'').

Our goal is to produce a Step~1 object $(X,\omega,Z,V^+,f,p)$ related to
the given graph $G$ as follows:
\begin{itemize}
\item $Z = F_1 \cup \ldots \cup F_n$, where each surface $F_i$
corresponds to a vertex $v_i$ in $G$ and the intersections correspond
to the edges in $G$.
\item $\genus(F_i) = g_i$
\item $\partial F_i$ has $m_i + d_i$ components.
\item For each $F_i$, there is one component $K$ of $\partial
F_i$ for which $\pf(K) = \cf(K) - d_i$, and for all other
components the two framings are equal.
\end{itemize}

Topologically $X$ is built from $0$--handles and $1$--handles, with one
$0$--handle for each vertex and each edge in $G$, and with a $1$--handle
connecting an edge $0$--handle to a vertex $0$--handle if that edge is
incident with that vertex.

We begin with two basic Step~1 objects $A =
(X_A,\omega_A,Z_A,V_A^+,f_A,p_A)$ and $B =
(X_B,\omega_B,Z_B,V_B^+,f_B,p_B)$ defined as follows (here we use
polar coordinates $(r_1,\theta_1,r_2,\theta_2)$ on $\reals^4$):
\begin{itemize}
\item $X_A = X_B = B^4 = \{r_1^2 + r_2^2 \leq 1\} \subset \reals^4$.
\item $\omega_A = \omega_B = r_1 dr_1 d\theta_1 + r_2 dr_2 d\theta_2$.
\item $Z_A = \{r_2 = 0\}$ and $Z_B = \{r_2 = 0\} \cup \{r_1 = 0\}$.
\item $V_A^+ = V_B^+ =\frac{1}{2}(r_1 \partial_{r_1} + r_2
\partial_{r_2})$.
\item $f_A = f_B = r_1^2 + r_2^2$.
\item $p_A = \theta_2$ while $p_B = \theta_1 + \theta_2$.
\end{itemize}
Thus $(\partial Z_A, p_A)$ is the standard open book on $S^3$ with
page equal to a disk and binding the unknot, and $(\partial Z_B, p_B)$
is the open book on $S^3$ with page equal to an annulus (a
left-twisted Hopf band), monodromy equal to a single right twist about
the core circle of the annulus, and the Hopf link with positive
linking number as the binding. Note that, for the single component of
$\partial Z_A$, we have $\pf = \cf$, whereas for each of the two
components of $\partial Z_B$, we have $\pf = \cf - 1$.

We think of these two objects as $4$--dimensional symplectic
$0$--handles, in the sense of Weinstein~\cite{Weinstein}. We also have
Weinstein's $4$--dimensional symplectic $1$--handle, which is
constructed as a neighborhood of the origin in $\reals^4$ with the
standard symplectic form $\omega = dx_1 dy_1 + dx_2 dy_2$, the Morse
function $f = -x_1^2 + y_1^2 + x_2^2 + y_2^2$ and the symplectic
dilation $V^+ = -x_1 \partial_{x_1} + 2 y_1 \partial_{y_1} +
\frac{1}{2}(x_2 \partial_{x_2} + y_2 \partial_{y_2})$. Weinstein shows
that we can always attach such a $1$--handle at any two points on a
convex boundary of a symplectic $4$--manifold, such that the symplectic
forms and symplectic dilations match up along the
glueing. In~\cite{GayCF} we have shown that, if in addition the
contact structure on the boundary is supported by an open book and the
attaching $3$--balls of the $1$--handle are pierced by the binding (as
the $z$--axis pierces the unit ball in $\reals^3$), then the handle can
be constructed in such a way that the new contact structure produced
by the associated contact surgery is also supported by an open
book. The new page is produced from the old page by attaching a
$2$--dimensional $1$--handle at the corresponding intervals along the
binding and the new monodromy is equal to the old monodromy extended
by the identity on the $1$--handle. (Note that the $2$--dimensional
$1$--handle is explicitly the set $\{x_2 = y_2 =0\}$ inside the
$4$--dimensional $1$--handle, and that this a symplectic surface.)

Now suppose that we are attaching such a $1$--handle to a Step~1 object
\linebreak
$(X,\omega,Z,V^+,f,p)$ at two points along the binding $\partial Z$
(not connecting two surfaces in $Z$ that already intersect). Then we
produce a new Step~1 object $(X_1,\omega_1,Z_1,V^+_1,f_1,p_1)$, where
$Z_1$ is the result of attaching a $1$--handle to $Z$, $(\partial
Z_1,p_1)$ is the open book described in the preceding paragraph, and
$f_1$ has a single new index $1$ critical point.

\begin{figure} 
\begin{center}
\includegraphics[width=4.5in,height=2in]{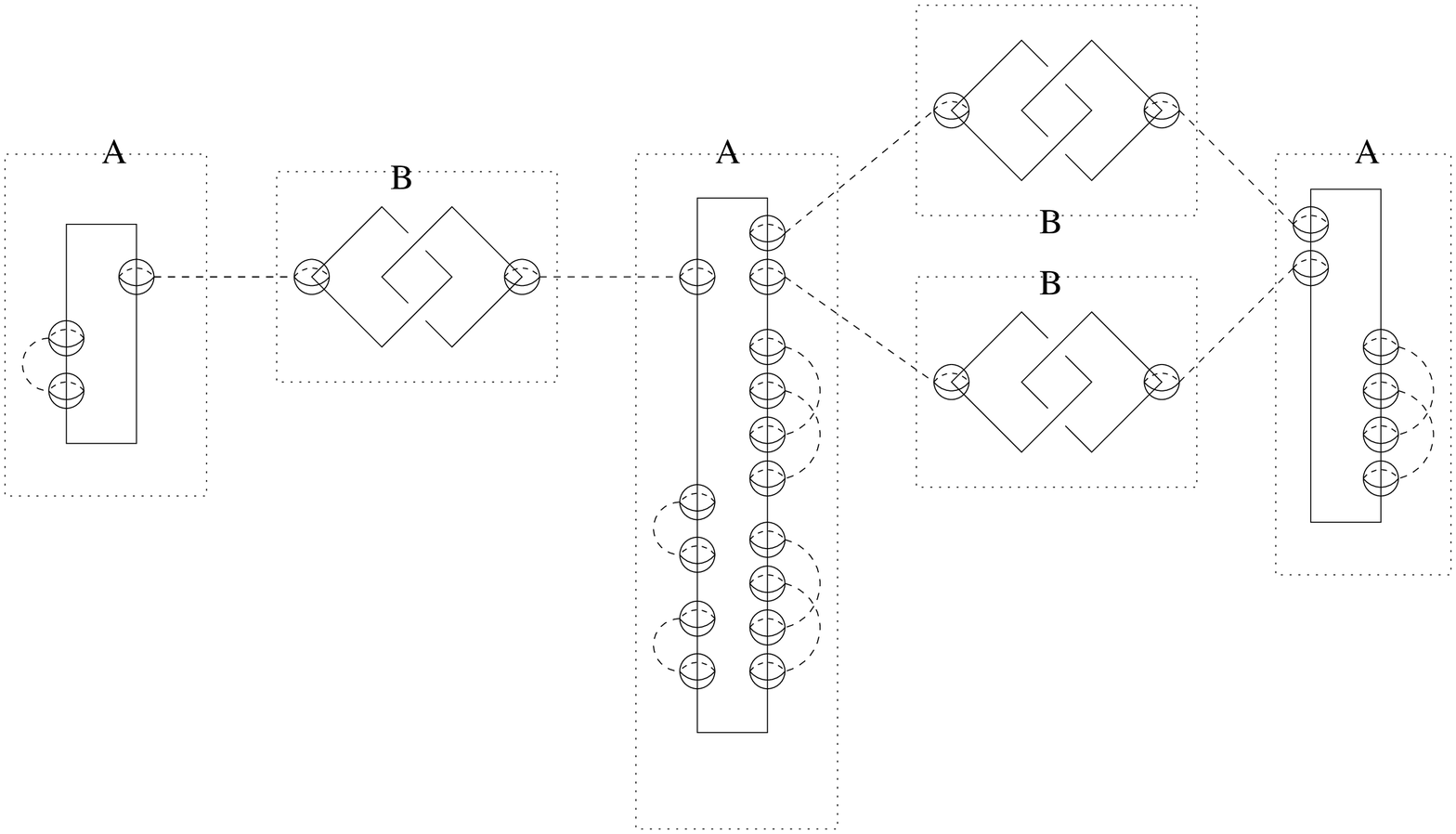}
\caption{The Step~1 object corresponding to the graph $G$ in
figure~\ref{F:SigmaGExample}}
\label{F:ManyS3s}
\end{center}
\end{figure}
We now describe how to build our desired Step~1 object corresponding
to the given graph $G$. Figure~\ref{F:ManyS3s} illustrates this
construction applied to the graph $G$ of figure~\ref{F:SigmaGExample}.
Start with a Step~1 object which is a disjoint union of many copies of
the $0$--handle objects $A$ and $B$, with one $A$ for each vertex and
one $B$ for each edge in $G$. These are indicated in
figure~\ref{F:ManyS3s} by dotted rectangles, with the bindings of the
open books indicated as solid links (unknots or Hopf links).  For an
edge $e_k$ connecting two vertices $v_i$ and $v_j$, let $B_k$ be the
corresponding copy of $B$ and $A_i$ and $A_j$ be the corresponding
copies of $A$. Connect $B_k$ to $A_i$ by a symplectic $1$--handle as
described above, with the $1$--handle connecting one component of the
binding in $B_k$ to the binding in $A_i$, and connect $B_k$ to $A_j$
by a $1$--handle connecting the other component of the binding in $B_k$
to the binding in $A_j$. In figure~\ref{F:ManyS3s} we have drawn the
attaching $3$--balls of these $1$--handles with dashed lines connecting
them. After doing this for all the edges, we have a Step~1 object
where the configuration $Z$ is a configuration of disks, one disk for
each vertex, with intersections given by $G$. Furthermore, for each
component $K_i$ of $\partial Z$ corresponding to a vertex $v_i$, we
have $\pf(K_i) = \cf(K_i) - d_i$. Now attach $2g_i$ $1$--handles in
pairs along each component $K_i$ of $\partial Z$ to get a Step~1
object for which the configuration is now a configuration of surfaces
$Z = F_1 \cup \ldots \cup F_n$ where each $F_i$ corresponds to a
vertex $v_i$, each $F_i$ has a single boundary component, and
$\genus(F_i) = g_i$, with the intersections given by $G$. These pairs
of $1$--handles are drawn on the lower right of each $A$ binding. We
still have $\pf(\partial F_i) = \cf(\partial F_i) - d_i$. Finally,
attach $(m_i + d_i - 1)$ $1$--handles along an isolated stretch of
$\partial F_i$ to get $m_i + d_i$ binding components for each $F_i$;
for the ``new'' binding components we will have $\pf = \cf$, while one
binding component still has $\pf = \cf - d_i$. These $1$--handles are
drawn on the lower left of each $A$ binding.
\begin{figure} 
\begin{center}
\includegraphics[width=4.5in, height=2in]{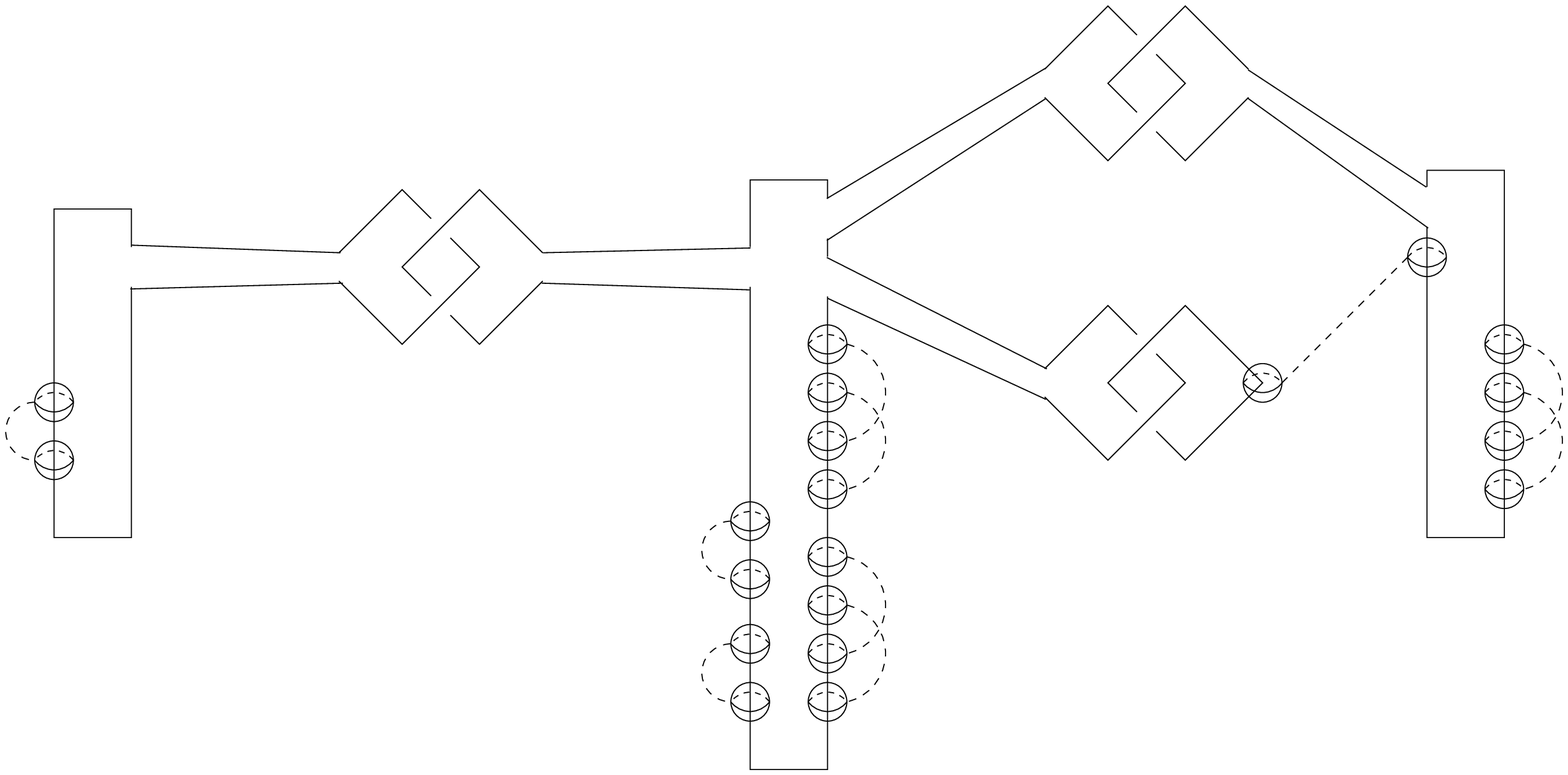}
\caption{The same Step~1 object, drawn as a Kirby calculus diagram in
a single $S^3$}
\label{F:OneS3}
\end{center}
\end{figure}
Figure~\ref{F:OneS3} shows a standard Kirby calculus diagram for the
same construction, drawn in a single copy of $S^3$; the link drawn is
the binding of an open book, not a surgery link (yet).

Note that the page of the resulting open book $(\partial C, p)$ in the
final Step~1 object $(X,\omega,Z,V^+,f,p)$ is exactly the surface
$\Sigma(G)$ associated to the graph $G$ and that the monodromy is
exactly the mapping class $\sigma(G)$. Thus $(\partial X, \xi^+) =
\mathcal{B}(\Sigma(G),\sigma(G))$.

{\bf Step 2}\qua For each $F_i$ in $Z = F_1 \cup \ldots \cup F_n$, choose
a positive constant $b_i < a_i/(2 \pi (m_i + d_i))$. By an explicit
construction (similar to that in section~4 of~\cite{GayCF}) one
can show that there exists a positive contact form $\alpha^+$
supported by $(\partial Z,p)$ with the following behavior near each
component $K$ of each $F_i$: In a neighborhood of $K$ there should
exist solid torus coordinates $(r,\mu,\lambda)$ (radial, meridinal and
longitudinal coordinates, with $K = \{r=0\})$ with respect to which
$\alpha^+ = \half r^2 (d\mu - d\lambda) + b_i d\lambda$ and $p = \mu +
\lambda$. After an isotopy fixing $\partial Z$, by theorem~\ref{TWTG},
we can assume that $\ker \alpha^+ = \xi^+ = \ker(\imath_{(V^+)}
\omega|_{\partial X})$.  Because $V^+$ is defined on all of $X$ and is
gradient-like for $f$, we can use the canonical symplectification of
$\xi^+$ and enlarge and/or trim $X$ so that in fact $\alpha^+ =
\imath_{(V^+)} \omega|_{\partial X}$.

At this point note that the area of each $F_i$ is $\int_{F_i} \omega =
\int_{\partial F_i} \alpha^+ = 2\pi (m_i + d_i) b_i < a_i$.

Now we recall some definitions from~\cite{Gay2Handles}. A {\em contact
pair} on a $3$--manifold $M$ is a pair $(\alpha^+,\alpha^-)$ of
$1$--forms defined, respectively, on open subsets $M^\pm$ with $M = M^+
\cup M^-$, such that $\pm \alpha^\pm \wedge d\alpha^\pm > 0$ on
$M^\pm$ and such that $d\alpha^+ = - d\alpha^-$ on $M^0 = M^+ \cap
M^-$. In particular $\alpha^+$ (resp.\ $\alpha^-$) is a positive
(resp.\ negative) contact form and $\alpha^0 = \alpha^+ + \alpha^-$ is
a closed, nowhere zero $1$--form on $M^0 = M^+ \cap M^-$. A {\em
dilation-contraction pair} on a symplectic $4$--manifold $(X,\omega)$
is a pair $(V^+,V^-)$ defined, respectively, on open subsets $X^\pm
\subset X$, such that $\mathcal{L}_{(V^\pm)} \omega = \pm \omega$ and
$\omega(V^+,V^-) = 0$. We say that $(V^+,V^-)$ {\em transversely
covers} a $3$--dimensional submanifold $M$ if $M \subset X^+ \cup X^-$
and both $V^+$ and $V^-$ are positively transverse to $M$. This gives
us an induced contact pair $(\alpha^+,\alpha^-)$ on $M$ defined by
$\alpha^\pm = \imath_{(V^\pm)} \omega|_M$, with domains $M^\pm = M \cap
X^\pm$.

In our situation we have the positive contact form $\alpha^+$ defined
on $(\partial X)^+ = \partial X$, supported by $(\partial Z,p)$; for a
large enough constant $k>0$, $(\alpha^+,\alpha^- = k dp - \alpha^+)$
will be a contact pair (with $(\partial X)^- = \partial X \setminus
\partial Z$). Fix such a $k$ and henceforth let $\alpha^- = k dp -
\alpha^+$.

Backward flow along the symplectic dilation $V^+$ starting on
$\partial X$ gives us an embedding $\phi^+: (-\infty,0] \times
\partial X \into X$ such that $\phi^+(0,p) = p$, $D\phi^+(\partial_t)
= V^+$ and $(\phi^+)^* \omega = d(e^t \alpha^+)$, where $t$ is the
coordinate on $(-\infty,0]$. Since $V^+$ is tangent to $Z$ and
gradient-like for $f$, we know that $\phi^+((-\infty,0] \times
(\partial X \setminus \partial Z)) = X \setminus Z$.  Lemma~4.1
in~\cite{Gay2Handles} then tells us that there exists a unique
symplectic contraction $V^-$ defined on $X \setminus Z$ such that
$(V^+,V^-)$ is a dilation-contraction pair transversely covering
$\partial X$ inducing the contact pair $(\alpha^+,\alpha^-)$. Forward
flow along $-V^-$ starting on $\partial X \setminus \partial Z$ then
gives an embedding $\phi^-$ from $\{(t,p) | 0 \leq t < F(p) \} \subset
\reals \times (-(\partial X \setminus \partial Z))$ into $X$ for some
function $F: \partial X \setminus \partial Z \to (0,\infty]$, such
that $\phi^-(0,p) = p$, $D\phi^-(\partial_t) = -V^-$ and $(\phi^-)^*
\omega = d(e^t \alpha^-)$. The proof of lemma~4.1
in~\cite{Gay2Handles} shows us how to explicitly calculate $V^-$ given
$(\alpha^+,\alpha^-)$, from which we can see that, in our case, $F <
\infty$ and the image of $\phi^-$ is all of $X \setminus Z$, with $Z
\setminus \partial Z = \{\lim_{t \to F(p)} \phi^-(t,p) | p \in
\partial X \setminus \partial Z \}$.

{\bf Step 3}\qua Our final symplectic $4$--manifold $(N(G),\omega(G))$
will be constructed by attaching a symplectic $2$--handle as described
in~\cite{Gay2Handles} along each component $K$ of the binding
$\partial Z \subset \partial X$ with framing $\pf(K)+1$, and then
attaching an open symplectic collar on the
boundary. In~\cite{Gay2Handles} it is shown that one can attach
handles in this way to produce a concave boundary, and
in~\cite{GayCF}, it is shown that the contact manifold on the boundary
is indeed $\mathcal{B}(\Sigma(G),\sigma(G)^{-1} \circ \delta(G))$. It
is not hard to see that the $4$--manifold produced in this way is a
neighborhood of a configuration of closed surfaces with the correct
genera and intersecting according to the graph $G$; the surfaces are
simply the surfaces $F_i$ in $Z$ capped off with the descending disks
of the $4$--dimensional $2$--handles.

To see that the self-intersections are correct, note that the
self-intersection of a surface $\Sigma_i$ built by attaching
$2$--handles, framed as above, along $\partial F_i$ for one of the
surfaces $F_i \subset Z$ is equal to the sum over all components $K$
of $\partial F_i$ of $(\pf(K)+1) - \cf(K)$, which is exactly $m_i$.

To see that the remaining claims of the theorem hold, we need to look
more closely at the structure of the $2$--handles.  Proposition~4.6
in~\cite{Gay2Handles} shows us how to construct our $2$--handles; here
we give the construction tailored to the special case at hand. For
each $F_i \subset Z$, let $c_i = a_i/(2\pi k (m_i+d_i))$ (with $k$ as
fixed in Step~2). The handle $H$ that will attach to each component
$K$ of $\partial F_i$ is a subset of $\reals^4$ with the symplectic
form $\omega_H = c_i (r_1 dr_1 d\theta_1 + r_2 dr_2 d\theta_2)$ with
the Morse function $f_H = -r_1^2 + r_2^2$. The following is a
dilation-contraction pair on $(\reals^4,\omega_H)$:
\[ V^+_H = (\half r_1 - \frac{k}{r_1}) \partial_{r_1} + \half r_2
\partial_{r_2} \]
\[ V^-_H = -\half r_1 \partial_{r_1} - (\half r_2 - \frac{k}{r_2})
\partial_{r_2} \] 

$(V^+_H,V^-_H)$ transversely covers the regular level sets of $f$ as
long as $-2k < f < 2k$. Let $\epsilon_1 = 2(b_i/c_i - k)$; note that
$-2k < \epsilon_1 < 0$. Choose some $\epsilon_2$ with $0 < \epsilon_2
< 2k$. Construct $H$ so that the attaching boundary of $H$ (which we
call $\partial_1 H$) is a neighborhood of $\{r_2 = 0\}$ in
$f^{-1}(\epsilon_1)$ and so that the free boundary $\partial_2 H$
interpolates from a neighborhood of $\{r_1 = 0\}$ in
$f^{-1}(\epsilon_2)$ down to $f^{-1}(\epsilon_1)$, so that both
boundaries are transverse to both $V_H^+$ and $V_H^-$ (where
defined). (See section~2 of~\cite{Gay2Handles} for a more detailed
discussion of this type of handle construction and notation.) On
$\partial_1 H$, we use solid torus coordinates $(r = \sqrt{c_i} r_2,
\mu = \theta_2, \lambda = -\theta_1)$; the contact pair induced by
$(V_H^+,V_H^-)$ on $\partial_1 H$ then becomes:
\[ (\alpha_H^+ = \half r^2 (d\mu - d\lambda) + b_i d\lambda,
\alpha_H^- = k(d\mu + d\lambda) - \alpha_H^+) \]
This is exactly the contact pair we have on a neighborhood of $K
\subset \partial X$, so that lemma~4.1 in~\cite{Gay2Handles} tells us
we can attach $H$ to $X$ by identifying the $(r,\mu,\lambda)$
coordinates on $\partial_1 H$ with the $(r,\mu,\lambda)$ coordinates
in a neighborhood of $K$, in such a way that the symplectic forms and
the dilation-contraction pairs fit together smoothly. Note that since
$p = \mu + \lambda$ in a neighborhood of $K$, we are attaching $H$
with framing $\pf(K) + 1$.

After attaching these handles to each component of $\partial F_i$, we
get a closed surface $\Sigma_i$ which is the union of $F_i$ and the
disks $D = \{r_2 = 0\} \cap H = \{r_2 = 0, r_1^2 \leq 2(k -
b_i/c_i)\}$ in each handle $H$. $\Sigma_i$ is smooth because $V^+$ is
tangent to $F_i$ and $V^+_H$ is tangent to $D$. $\Sigma_i$ is
symplectic because $F_i$ and $D$ are both symplectic.  We have already
arranged that the area of $F_i$ is $2 \pi (m_i+d_i) b_i$. The area of
each disk $D$ is $2 \pi c_i (k - b_i/c_i)$. Thus the area of
$\Sigma_i$ is exactly $a_i$. Let $C(G) = \Sigma_1 \cup
\ldots \cup \Sigma_n$.

Note that, in the handle $H$, the symplectic contraction $V^-_H$ is
defined across all of the free boundary $\partial_2 H$ whereas the
dilation $V^+_H$ does not extend across the ascending circle $\{r_1 =
0\} \cap \partial_2 H$. Thus after attaching all the handles we get a
symplectic $4$--manifold $(X_1,\omega_1)$ with a dilation-contraction
pair $(V^+_1,V^-_1)$ which transversely covers $\partial X_1$ inducing
a contact pair $(\alpha_1^+,\alpha_1^-)$ with domains $(\partial
X_1)^+ =\partial X_1 \setminus L_1$ (where $L_1$ is the union of the
ascending circles) and $(\partial X_1)^- = \partial X_1$. The closed
$1$--form $\alpha_1^0 = \alpha_1^+ + \alpha_1^-$ is $k dp_1$ for the
natural open book $(L_1,p_1)$ on $\partial X_1$ that results from
$\pf+1$ surgeries on the binding of the open book $(\partial Z,p)$ on
$\partial X$, the pages of which are still diffeomorphic to
$\Sigma(G)$ and the monodromy of which is now $\sigma(G) \circ
\delta(G)^{-1}$ (see~\cite{GayCF}). The fact that
$(\alpha_1^+,\alpha_1^-)$ is a contact pair implies that $\alpha_1^0
\wedge \alpha_1^- > 0$, which in turn implies that $\alpha_1^-$, as a
positive contact form on $-\partial X_1$, is supported by the open
book $(L_1,-p_1)$, which has page $\Sigma(G)$ and monodromy
$\sigma(G)^{-1} \circ \delta(G)$. Thus we let $M(G) = - \partial X_1$
and $\xi(G) = \ker \alpha^-$. Also let $\alpha = \alpha^-$.

From the explicit form for $V^-_H$ and the observations made in Step~2
about flow along $-V^-$, we see that flow along $-V^-_1$ starting on
$\partial X_1$ gives a diffeomorphism $\phi^-_1$ from $\{(t,p) | 0
\leq t < f(p)\} \subset \reals \times M(G)$ to $X_1 \setminus C(G)$,
for some function $f: M(G) \to (0,\infty)$, such that $\phi^-(0,p) =
p$, $D\phi^-(\partial_t) = -V^-$, $(\phi^-)^* \omega_1 = d(e^t
\alpha)$ and $C(G) = \{ \lim_{t \to f(p)} \phi^-(t,p) | p \in M(G)
\}$.

Finally let $(N(G),\omega(G)) = (X_1,\omega_1) \cup_{\phi^-}
((-\infty,0] \times M(G), d(e^t \alpha))$; the embedding $\phi$ is
simply $\phi^-$ extended by the identity on $(-\infty,0] \times M(G)$.

\end{proof}

\section{The nonpositive case}

If our initial graph $G$ is not positive, we can simply add extra
vertices labelled $(0,0,1)$ to produce a graph $G' \supset G$ which is
positive. This corresponds to plumbing on some extra spheres of square
$0$. We can now carry out the construction above applied to $G'$, but
stop short of attaching the $2$--handles required to close off these
extra spheres. This will give a model neighborhood $(N(G),\omega(G))$
of a configuration $C(G)$ with $G = G(C(G))$, but now the boundary
will not be concave. Instead the boundary will be ``partially convex
and partially concave'' in the following sense: $(N(G),\omega(G))$
will carry a dilation-contraction pair transversely convering
$\partial N(G)$ inducing a contact pair $(\alpha^+,\alpha^-)$, but
neither $\alpha^+$ nor $\alpha^-$ will be defined on all of $\partial
N(G)$. Nevertheless, the pair $(\alpha^+,\alpha^-)$ will determine the
germ of $\omega(G)$ along $\partial N(G)$ and the dilation-contraction
pair determines something like a canonical symplectification of the
contact pair, so that we have good control on the symplectic topology
of $N(G) \setminus C(G)$. Furthermore, $(\alpha^+,\alpha^-)$ will be
supported by a {\em signed open book} $(L^+,L^-,p)$ on $M = \partial
N(G)$, by which we mean the following: $(L = L^+ \amalg L^-,p)$ is an
open book, $\alpha^\pm$ is defined on $M^\pm = M \setminus L^\mp$, and
$\alpha^\pm$ is supported, as a positive contact form on $\pm M^\pm$,
by $(L^\pm, \pm p)$. The link $L^-$ will be the union of the ascending
circles for the $2$--handles that we did attach, while $L^+$ will be
the binding components from the Step~1 object to which we did not
attach $2$--handles. Smaller or larger neighborhoods of $C(G)$ will
have boundaries which are still transversely covered by the
dilation-contraction pair and hence carry related contact pairs, all
supported by the same signed open book. This line of reasoning will be
investigated more thoroughly in a future paper.

\section{Boundaries of arbitrary plumbings}

\begin{proof}[Proof of theorem~\ref{TopThm}]
If we strip the symplectic topology out from the proof of
theorem~\ref{MainThm}, Step~2 is irrelevant. In Step~1, relabel $B$ as
$B_+$ and introduce a negative version of $B$ which we call $B_-$;
$B_+$ (resp.\ $B_-$) is a neighborhood of a positive (resp.\ negative)
intersection of two disks, with a left-twisted (resp.\ right-twisted)
Hopf band as the open book on the boundary, with monodromy equal to a
single right (resp.\ left) Dehn twist along the core of the band. On
$\partial B_+$ we have $\pf = \cf -1$ and on $\partial B_-$ we have
$\pf = \cf + 1$. Thus if we mimic the construction in
theorem~\ref{MainThm} but use copies of $B_+$ for $+$ edges and copies
of $B_-$ for $-$ edges, we can produce a (non-symplectic) Step~1
object $(X,Z,f,p)$ where:
\begin{itemize}
\item $X$ is a compact $4$--manifold with boundary which is a
neighborhood of $Z$.
\item $Z = F_1 \cup \ldots \cup F_n \subset X$ is a configuration of properly
embedded surfaces with boundary, corresponding to the vertices of $G$,
with the appropriate genera and self-intersections and intersecting
transversely according to the edges of $G$.
\item $\partial F_i$ has $|m_i + d_i|$ components, unless $m_i + d_i =
0$, in which case $\partial F_i$ has one component. (Recall that now
$d_i$ is the signed count of edges connecting to $v_i$.)
\item $p : \partial X \setminus \partial Z \rightarrow S^1$ is a
fibration making $(\partial Z,p)$ an open book on $\partial X$.
\item For each $F_i$, there is one component $K$ of $\partial F_i$ for
which $\pf(K) = \cf(K) - d_i$ and for all other components the two
framings are equal.
\end{itemize}
When we get to Step~3, since we are no longer requiring that our
handles be symplectic, we can attach $2$--handles along binding
components with any framings we choose. Framing $\pf-1$ produces a new
open book with the same page and introduces a right-handed boundary
Dehn twist into the monodromy; framing $\pf + 1$ also produces a new
open book with the same page and introduces a left-handed boundary
Dehn twist. Framing $\pf$ produces a new open book with the page
alterred by capping off the corresponding boundary component,
decreasing the number of binding components by $1$. For vertices $v_i$
with $m_i + d_i < 0$, use $(\pf -1)$-framed binding handles, for
vertices with $m_i + d_i = 0$, use $\pf$-framed binding handles and
for vertices with $m_i + d_i > 0$, use $(\pf + 1)$-framed binding
handles. Because of these choices of framings we then calculate that
for each $\Sigma_i$, $\Sigma_i \cdot \Sigma_i = m_i$.
\end{proof}

%
%
%

\Addresses\recd

\newpage


\setcounter{section}{0}

\lognumber{45}
\volumenumber{3}
\volumeyear{2003}
\papernumber{45}
\published{20 December 2003}
\pagenumbers{1275}{1276}
\received{1 December  2003}

\title{Correction to ``Open books and\\configurations of symplectic
  surfaces''}                    

\def\theaddress{CIRGET, Universit\'{e} du Qu\'{e}bec \`{a} Montr\'{e}al\\
Case Postale 8888, Succursale centre-ville\\ 
Montr\'{e}al (QC) H3C 3P8, Canada}                  

\asciiaddress{CIRGET, Universite du Quebec a Montreal\\
Case Postale 8888, Succursale centre-ville\\ 
Montreal (QC) H3C 3P8, Canada}

\def\theemail{gay@math.uqam.ca}

\begin{abstract} 
We correct the main theorem in ``Open books and configurations of
symplectic surfaces''~\cite{gay} and its proof. As originally stated,
the theorem gave conditions on a configuration of symplectic surfaces
in a symplectic $4$-manifold under which we could construct a model
neighborhood with concave boundary and describe explicitly the open
book supporting the contact structure on the boundary. The statement
should have included constraints on the areas of the surfaces.
\end{abstract}

\primaryclass{57R17}                
\secondaryclass{57N10, 57N13}              
\maketitle

In the paper being corrected~\cite{gay}, we considered symplectic
configuration graphs where each vertex was decorated with a triple
$(g_i,m_i,a_i)$, $g_i$ being the genus of a surface, $m_i$ its
self-intersection, and $a_i$ its symplectic area. Theorem~1.1
in~\cite{gay} is false as stated. However we have the following:
\begin{correction}
  The conclusions of theorem~1.1 (parts~A and ~B) are true if we add
  the hypothesis that there exists a constant $\lambda > 0$ such that,
  for each $i$, $a_i = \lambda(m_i + d_i)$.
\end{correction}

An easy counterexample to the theorem as originally stated is given by
two surfaces $\Sigma_1$ and $\Sigma_2$, with $\Sigma_1 \cdot \Sigma_1
= \Sigma_2 \cdot \Sigma_2 = 1$ and $\Sigma_1 \cdot \Sigma_2 = 1$. If
$\int_{\Sigma_1} \omega = a_1 \neq \int_{\Sigma_2} \omega = a_2$, then
$[\omega] \neq 0$ when restricted to the boundary of a neighborhood of
$\Sigma_1 \cup \Sigma_2$, because $[\omega]$ is nonzero when evaluated on
the class $[\Sigma_1] - [\Sigma_2]$, which lives on the boundary. In
general, this problem occurs when the intersection matrix for the
configuration of surfaces has determinant equal to zero.

The error is on line~4 of page~583 in~\cite{gay}, in the calculation
of $\alpha_H^-$. The correct statement is:
\[ \alpha_H^- = k c_i (d \mu + d\lambda) - \alpha_H^+ \]
Thus, to arrange that this agree with the contact pair we already have
on a neighborhood of $K \subset \partial X$, we are forced to choose
$c_i = 1$, as opposed to $c_i = a_i/(2 \pi K(m_i+d_i))$, which was the
choice made on page~582 line~22. With $c_i = 1$, we get that the area
of $\Sigma_i$ is $2\pi (m_i + d_i) k$. After the construction we can
rescale by a constant to get the area to be $\lambda (m_i + d_i)$.

Note that, if the configuration $\Sigma_1 \cup \ldots \cup \Sigma_n$
is contained in a closed symplectic $4$-manifold $(X,\omega)$, then
the area condition added in this correction is equivalent to the
condition that $[\omega]$ is Poincar\'{e} dual to some multiple of
$[\Sigma_1] + \ldots + [\Sigma_n] + \gamma$, where $\gamma \in
H_2(X;\mathbb{Z})$ is a class with $\gamma \cdot \Sigma_i = 0$ for $i
= 1, \ldots, n$. This is because $m_i + d_i = \Sigma_1 \cdot \Sigma_i
+ \ldots + \Sigma_n \cdot \Sigma_i$.

Here we emphasize that this correction does not affect the
applications in section~2 of~\cite{gay}, or forthcoming applications
in~\cite{gaykirby}.

\Addresses\recd
\end{document}


\begin{thebibliography}

\bibitem{AkbOzb}
\textbf{Selman Akbulut}, \textbf{Burak Ozbagci}, \emph{On the topology of
  compact {S}tein surfaces}, Int. Math. Res. Not.  (2002) 769--782

\bibitem{EliashFilling}
\textbf{Yakov Eliashberg}, \emph{Filling by holomorphic discs and its
  applications}, from: ``Geometry of low-dimensional manifolds, 2 (Durham,
  1989)'', Cambridge Univ. Press, Cambridge (1990)  45--67

\bibitem{GayCF}
\textbf{David~T Gay}, \emph{Explicit concave fillings of contact
  three-manifolds}, Math. Proc. Cambridge Philos. Soc. 133 (2002) 431--441

\bibitem{Gay2Handles}
\textbf{David~T Gay}, \emph{Symplectic 2-handles and transverse links}, Trans.
  Amer. Math. Soc. 354 (2002) 1027--1047 (electronic)

\bibitem{Giroux}
\textbf{Emmanuel Giroux}, in preparation

\bibitem{Gromov}
\textbf{M Gromov}, \emph{Pseudoholomorphic curves in symplectic manifolds},
  Invent. Math. 82 (1985) 307--347

\bibitem{LoiPier}
\textbf{Andrea Loi}, \textbf{Riccardo Piergallini}, \emph{Compact {S}tein
  surfaces with boundary as branched covers of ${B}\sp 4$}, Invent. Math. 143
  (2001) 325--348

\bibitem{McDuffCvx}
\textbf{Dusa McDuff}, \emph{Symplectic manifolds with contact type boundaries},
  Invent. Math. 103 (1991) 651--671

\bibitem{McDSal}
\textbf{Dusa McDuff}, \textbf{Dietmar Salamon}, \emph{Introduction to
  symplectic topology}, second edition, Oxford University Press, Oxford (1998)

\bibitem{Symington}
\textbf{Margaret Symington}, \emph{Symplectic rational blowdowns}, J.
  Differential Geom. 50 (1998) 505--518

\bibitem{SymingtonGen}
\textbf{Margaret Symington}, \emph{Generalized symplectic rational blowdowns},
  Algebr. Geom. Topol. 1 (2001) 503--518 (electronic)

\bibitem{ThurstonWinkel}
\textbf{W\,P Thurston}, \textbf{H\,E Winkelnkemper}, \emph{On the existence of
  contact forms}, Proc. Amer. Math. Soc. 52 (1975) 345--347

\bibitem{Torisu}
\textbf{Ichiro Torisu}, \emph{Convex contact structures and fibered links in
  3-manifolds}, Internat. Math. Res. Notices 2000 (2000) 441--454

\bibitem{Wajnryb}
\textbf{Bronislaw Wajnryb}, \emph{An elementary approach to the mapping class
  group of a surface}, Geom. Topol. 3 (1999) 405--466 (electronic)

\bibitem{Weinstein}
\textbf{Alan Weinstein}, \emph{Contact surgery and symplectic handlebodies},
  Hokkaido Math. J. 20 (1991) 241--251

\end{thebibliography}

\begin{thebibliography}


\bibitem{gay}
{\bf D T Gay}
{\em  Open books and configurations of symplectic surfaces},
Alg. Geom. Top. 3 (2003), 569--586


\bibitem{gaykirby}
{\bf D T Gay and R Kirby}
{\em Constructing symplectic forms on 4-manifolds which vanish on
  circles},
to appear


\end{thebibliography}
\end{document}